\input dd.def
\documentclass{amsart}
\usepackage{amsmath}
\usepackage{amsfonts}
\usepackage{amssymb}
\usepackage{graphicx}
\usepackage{epsfig}
\newtheorem{theorem}{Theorem}[section]
\newtheorem{lemma}[theorem]{Lemma}

\theoremstyle{definition}

\theoremstyle{remark}
\newtheorem{remark}[theorem]{Remark}

\numberwithin{equation}{section}

\def\conj{\mathop{\fam0 conj}}
\def\Fix{\mathop{\fam0 Fix}}
\def\Vol{\mathop{\fam0 Vol}}
\def\Re{\mathop{\fam0 Re}}
\def\Im{\mathop{\fam0 Im}}

\makeatletter
\ifx\undefined\setboxz@h\gdef\setboxz@h{\setbox\z@\hbox}\fi
\ifx\undefined\boxz@\gdef\boxz@{\box\z@}\fi
\ifx\undefined\wdz@\gdef\wdz@{\wd\z@}\fi \makeatother

\def\bA{\mathbf{A}}
\def\bB{\mathbf{B}}
\def\er#1{E_{\mathbb{R}}^{(#1)}}
\def\BAR#1{\bar{#1}}
\def\+{\mathbin{\scriptstyle\sqcup}}
\def\eightpoint{\small}

\def\CP{\mathcal{P}}

\begin{document}

\title{Monodromy groups of Real Enriques Surfaces}
\author{Sultan Erdo\u{g}an Dem\.ir}
\address{bilkent university, department of mathematics, 06800 ankara, turkey}
\email{erdogan@fen.bilkent.edu.tr}
\thanks{The author is partially supported by T\"UB\.ITAK, The Scientific and Technological Research Council of Turkey}




\begin{abstract}
We compute the monodromy groups of real Enriques surfaces of
hyperbolic type. The principal tools are the deformation
classification of such surfaces and a modified version of
Donaldson's trick, relating real Enriques surfaces and real
rational surfaces.
\end{abstract}

\maketitle
\section{Introduction}
An {\em Enriques surface} is a complex analytic surface $E$ with
$\pi_1(E)=\mathbb{Z}_2$ and having a $K3$-surface $X$ as its
universal cover. An Enriques surface is called {\em real} if it is
supplied with an anti-holomorphic involution $\conj$, called {\em
complex conjugation}. The {\em real part} of a real surface $E$ is
the fixed point set $E_{\mathbb{R}} =\Fix \conj$. A {\em
topological type} of real surfaces is a class of surfaces with
homeomorphic real parts. A real Enriques surface $E$ is a smooth
$4$-manifold, its real part $E_{\mathbb{R}}$ is either empty or a closed
$2$-manifold with finitely many components, each being either
$S=S^2$, or $S_g=\sharp_g(S^1 \times S^1)$, or $V_p= \sharp_p
\mathbb{RP}^2$.

Let $E$ be a real Enriques surface and $p:X\rightarrow E$ its
universal covering. Denote by $\tau : X\rightarrow X$ the deck
translation of $p$, called the {\em Enriques involution}. There
are exactly two liftings $t^{(1)},t^{(2)}: X \rightarrow X$ of
$\conj$ to $X$, which are both anti-holomorphic involutions. They
commute with each other and with $\tau$, and their composition is
$\tau$. For both $i=1,2$, the real parts
$X_{\mathbb{R}}^{(i)}=\Fix t^{(i)}$, and their images
$E_{\mathbb{R}}^{(i)}=p(X_{\mathbb{R}}^{(i)})$ (called the {\em
halves} of $E_{\mathbb{R}}$) are disjoint, $E_{\mathbb{R}}^{(i)}$
consists of whole components of $E_{\mathbb{R}}$, and
$E_{\mathbb{R}}=E_{\mathbb{R}}^{(1)} \sqcup E_{\mathbb{R}}^{(2)}$.
This decomposition is a deformation invariant of pair $(E,\conj)$.
We use the notation
$E_\mathbb{R}=\{\textrm{half}~E_{\mathbb{R}}^{(1)}\} \sqcup \{
\textrm{half}~E_{\mathbb{R}}^{(2)}\}$ for the half decomposition.
To describe the topological types of the real part the concept of
{\em topological Morse simplification, i.e.}, Morse transformation
of the topological type which decreases the total Betti number, is
used. A topological Morse simplification is either removing a
spherical component $(S\rightarrow \varnothing)$ or contracting a
handle $(S_{g+1}\rightarrow S_g \textrm{ or }V_{p+2}\rightarrow
V_p)$. The complex deformation type of surfaces being fixed
(\textit{e.g.}, $K3$ or Enriques), a topological type is called
{\em extremal} if it cannot be obtained from another one (in the
same complex deformation type) by a topological Morse
simplification.

The classification of real Enriques surfaces up to deformation was
given by A. Degtyarev, I. Itenberg and V. Kharlamov in~\cite{D2},
where one can find a complete list of deformation classes, the
invariants necessary to distinguish them, and detailed
explanations of the invariants. It turns out that the deformation
class of a real Enriques surface is determined by the topology of
its complex conjugation involution. Deformation classification can
be regarded as the study of the set of connected components, ({\em
i.e.}, $\pi_0$) of the moduli space. In this paper, we pose the
question about its fundamental group ({\em i.e.}, $\pi_1$). More
precisely, we study the canonical representation of the
fundamental group of a connected component of the moduli space in
the group~$G$ of permutations of the components of the real part
of the corresponding surfaces. In other words, we discuss the {\em
monodromy groups} of real Enriques surfaces, {\em i.e.}, the
subgroups of~$G$ realized by `auto-deformations' and/or
automorphisms of the surfaces.

The similar question for various families of $K3$-surfaces has been extensively
covered in the literature. Thus, the monodromy groups have been studied for
nonsingular plane sextics by Itenberg~\cite{I} and for nonsingular surfaces of
degree four in ${\mathbb{RP}}^3$ by Kharlamov~\cite{K1},~\cite{K2} and~\cite{K3} and Moriceau~\cite{M}.

A real Enriques surface is said to be of {\em hyperbolic,
parabolic,} or {\em elliptic type} if the minimal Euler
characteristic of the components of $E_{\mathbb{R}}$ is negative,
zero, or positive, respectively. In the deformation
classification, hyperbolic and parabolic cases are treated
geometrically (based on Donaldson's trick~\cite{Don}) whereas the
elliptic cases are treated arithmetically (calculations using the
global Torelli theorem for $K3$-surfaces~ {\em cf.}~\cite{BPV}).
There also is a crucial difference between the approaches to
surfaces of hyperbolic and parabolic types. In the former case,
natural \emph{complex} models of complex \textit{DPN}-pairs are
constructed, and a real structure descends to the model by
naturality. In the latter case, it is difficult to study complex
\textit{DPN}-pairs systematically and \emph{real} models of
\emph{real} \textit{DPN}-pairs are constructed from the very
beginning. We study the surfaces of hyperbolic types in this work.
Thus, we deal with an equivariant version of Donaldson's  trick
for Enriques surfaces modified by A. Degtyarev and V. Kharlamov~\cite{D2}, which transforms a real Enriques surface to a real
rational surface with a nonsingular real anti-bicanonical curve on
it. We analyze this construction and adopt it to the study of the
monodromy groups. In particular, we discuss the conditions
necessary for an additional automorphism of the rational surface
to define an automorphism of the resulting real Enriques surface.
The principal result of the paper can be roughly stated as follows
(for the exact statements see Theorems \ref{DP}, \ref{almostDP}
and \ref{(2,r)}): {\em For a real Enriques surface of hyperbolic
type, with some exceptions listed explicitly in each statement,
any permutation of homeomorphic components of each half of
$E_{\mathbb{R}}$ can be realized by deformations and/or
automorphisms}.

The exceptions deserve a separate discussion. In most cases, the
nonrealizable permutations are prohibited by a purely topological
invariant, the so-called Pontrjagin-Viro form (see~\cite{D4} and
remarks following the relevant statements). There are, however, a
few surfaces, those with $E_\mathbb{R}^{(1)}=V_3\sqcup ...$, for
which the Pontrjagin-Viro form is not well defined but the
spherical components of $E_\mathbb{R}^{(1)}$ cannot be permuted.
The question whether these permutations are realizable by
equivariant auto-homeomorphisms of the surface remains open. Same
question for parabolic and elliptic cases is a subject of a future
study as it seems to require completely different means.

Organization of the paper is as follows: In Section~\ref{Pre}, we
recall some rational surfaces, curves on them, and a few results
related to their classification up to rigid isotopy. In
Section~\ref{DPN}, we describe (modified) Donaldson's trick and
the resulting correspondence theorem, the construction backwards,
and recall some results concerning specific families of real
Enriques surfaces of hyperbolic type. In Section~\ref{Lifting}, a
few necessary conditions for lifting automorphisms are discussed.
In Section~\ref{Main}, the main result is stated and proved in
three theorems.

\textit{Acknowledgments.} This work is a part of the author's
PhD study. She would like to express her deepest gratitude to
her supervisor Prof. Alex Degtyarev for his guidance, patience and
invaluable advice.

\section{Some surfaces and curves on them}\label{Pre}

\subsection{\textit{DPN}-pairs}
A nonsingular algebraic surface admitting a nonempty nonsingular
anti-bicanonical curve ({\em i.e.,} curve in the class $|{-2K}|$),
is called a {\em DPN-surface}. Most \textit{DPN}-surfaces are
rational.

A pair $(Y,B)$, where $Y$ is a \textit{DPN}-surface and $B\in
|{-2K}|$ is a nonsingular curve, is called a {\em DPN-pair}. A
\textit{DPN}-pair $(Y,B)$ is called {\em unnodal} if $Y$ is
unnodal (does not contain a ${(-2)}$-curve), {\em rational} if $Y$
is rational, and {\em real} if both $Y$ and $B$ are real. The {\em
degree} of a rational \textit{DPN}-pair $(Y,B)$ is the degree of
$Y$, {\em i.e., $K^2$}.

If $(Y,B)$ is a rational \textit{DPN}-pair, the double covering
$X$ of $Y$ ramified along $B$ is a $K3$-surface. A
\textit{DPN}-surface contains finitely many $(-4)$-curves. A
rational \textit{DPN}-surface $Y$ of degree $d$ that has $r$
$(-4)$-curves is called a {\em $(g,r)$-surface}, where $g=d+r+1$.
In fact $g\geq 1$ and any nonsingular curve $B\in |{-2K_{Y}}|$ is
one of the following topological types (see~\cite{D2}):

\begin{enumerate}
\item $B\cong S_{g}\sqcup rS$ if $g>1$; \item $B\cong S_{1}\sqcup
rS$ or $rS$ if $g=1$ and $r>0$; \item $B\cong 2S_{1}$ or $S_{1}$
if $g=1$ and $r=0$.
\end{enumerate}

Let $Y$ be a real surface with $H_{1}(Y)=0$. An {\em admissible
branch curve} on $Y$ is a nonsingular real curve $B\subset Y$ such
that $[B]=0$ in $H_{2}(Y)$, the real part $B_{\mathbb{R}}$ is
empty and $B$ is not linked with $Y_{\mathbb{R}}$. An {\em
admissible DPN-pair} is a real rational {\em DPN}-pair $(Y,B)$
with $B$ an admissible branch curve.

Donaldson's trick (see subsection~\ref{Don}) establishes a
one-to-one correspondence between the set of deformation classes
of real Enriques surfaces with distinguished nonempty half ({\em
i.e.}, pairs $(E,E_{\mathbb{R}}^{(1)})$ with
$E_{\mathbb{R}}^{(1)}\neq \varnothing$) and the set of deformation
classes of admissible \textit{DPN}-pairs $(Y,B)$. Inverse
Donaldson's trick (see subsection~\ref{inverse}) establishes a
surjective map from the set of deformation classes of unnodal
admissible \textit{DPN}-pairs to the set of deformation classes of
real Enriques surfaces with distinguished nonempty half.

\subsection{Del Pezzo and geometrically ruled rational surfaces}
A {\em Del Pezzo surface} $Y$ is a surface such that ${K_Y^2>0}$
and ${D \cdot K_Y \leq 0}$ for any effective divisor $D$ on $Y$.
An {\em unnodal Del Pezzo surface} $Y$ is a surface whose
anticanonical divisor is ample, or equivalently, a Del Pezzo
surface without $(-2)$-curves.

We use the notation $\Sigma_a$, $a\geq 0$, for the {\em
geometrically ruled rational surface}  ({\em i.e.}, relatively
minimal conic bundle over $\mathbb{P}^1$) that has a section of
square $(-a)$, which is called the exceptional section. The
classes of the exceptional section $E_0$ and of a generic section
is denoted by $e_0$ and $e_{\infty }$, respectively, so that
$e_0^2={-a},$ $e_{\infty}^2=a,$ and $e_0\cdot e_{\infty}=0$. The
class of the fiber (generatrix) will be denoted by $l$; one has
$l^2=0$ and $l\cdot e_0= l\cdot e_{\infty}=1$. Any irreducible
curve in $\Sigma_a$ with $a\geq 1$, either is $E_0$ or belongs to
$|xl+ye_{\infty}|$, $x,y\geq 0$. If $a=0$ then $e_0= e_{\infty}$.
Thus, if $l_1$ denotes $e_0= e_{\infty}$ and $l_2$ denotes $l$
then any irreducible curve in $\Sigma_0$ belongs to $|xl_1+yl_2|$,
$x,y\geq 0$.

\subsection{Rigid isotopies} Recall that an {\em isotopy} is a homotopy from one
embedding of a manifold $M$ into a manifold $N$ to another
embedding such that, at every time, it is an embedding. An isotopy
in the class of nonsingular (or, more generally, equisingular, in
some appropriate sense) embeddings of analytic varieties is called
{\em rigid}. Below we are mainly dealing with rigid isotopies of
nonsingular curves on rational surfaces. Clearly, such an isotopy
is merely a path in the space of nonsingular curves.

An obvious rigid isotopy invariant of a real curve $C$ on a real
surface $Z$ is its {\em real scheme, i.e.,} the topological type
of the pair $(Z_{\mathbb{R}},C_{\mathbb{R}})$.

The deformation classification of real Enriques surfaces and hence
the monodromy problem of those leads to a variety of auxiliary
classification problems for curves on surfaces and surfaces in
projective spaces. Below we give a brief account of the related
results and recall the basic definitions and facts about them.
Details and further references can be found, {\em e.g.},
in~\cite{D2}.

\subsection{Curves in $\mathbb{P}_{\mathbb{R}}^2$}

The real point set $C_{\mathbb{R}}$ of a nonsingular curve $C$ in
$\mathbb{P}_{\mathbb{R}}^2$ is a collection of circles $A$
embedded in $\mathbb{P}_{\mathbb{R}}^2$, two- or one-sidedly. In
the former case the component is called an {\em oval}. Any oval
divides $\mathbb{P}_{\mathbb{R}}^2$ into two parts; the {\em
interior} of the oval, homeomorphic to a disk and the {\em
exterior} of the oval, homeomorphic to the M\"obius band. The
relation {\em to be in the interior of} defines a partial order on
the set of ovals, and the collection $A$ equipped with this
partial order determines the {\em real scheme} of $C$. The
following notation is used to describe real schemes: If a real
scheme has a single component, it is denoted by $\langle J
\rangle$, if the component is one-sided, or by $\langle 1
\rangle$, if it is an oval. The empty real scheme is denoted by
$\langle 0 \rangle$. If $\langle \mathcal{A} \rangle $ stands for
a collection of ovals, the collection obtained from it by adding a
new oval surrounding all the old ones is denoted by $\langle
1\langle \mathcal{A} \rangle \rangle$. If a real scheme splits
into two subschemes $\langle \mathcal{A}_1 \rangle $, $\langle
\mathcal{A}_2 \rangle $ so that no oval of $\langle \mathcal{A}_1
\rangle $ (respectively, $\langle \mathcal{A}_2 \rangle $)
surrounds an oval of $\langle \mathcal{A}_2 \rangle $
(respectively, $\langle \mathcal{A}_1 \rangle $), it is denoted by
$\langle \mathcal{A}_1\sqcup \mathcal{A}_2 \rangle $. If a real scheme
contains $n$ disjoint copies of $\langle 1 \rangle $ it is denoted by $\langle n \rangle $.

\begin{theorem} \label{p2}{\rm(\cite{K})}
A nonsingular real quartic $C$ in $\mathbb{P}^2$ is determined up
to rigid isotopy by its real scheme. There are six rigid isotopy
classes, with real schemes $\langle\alpha \rangle$,
$\alpha=0,...,4$ and $\langle 1 \langle 1 \rangle \rangle$.
\end{theorem}

\begin{lemma} \label{quartic}{\rm(\cite{D2})}
Let $C$ be a nonsingular real quartic with the real scheme
$\langle\alpha \rangle$, $\alpha=2,3,4$ in $\mathbb{P}^2$. Then
any permutation of the ovals of $C$ can be realized by a rigid
isotopy.
\end{lemma}

\subsection{Cubic sections on a quadratic cone} Let ${U\in
|ne_{\infty}|}$ be a nonsingular real curve in $\Sigma_2$ with its
standard real structure ($(\Sigma_2)_\mathbb{R}=S_1$). Each
connected component of $U_{\mathbb{R}}$ is either an oval or
homologous to $(E_0)_{\mathbb{R}}$. The latters, together with
$(E_0)_{\mathbb{R}}$, divide $(\Sigma_2)_{\mathbb{R}}$ into
several connected components $Z_1,...,Z_k$. Fixing an orientation
of the real part of a real generatrix of $\Sigma_2$ determines an
order of the components $Z_i$, and the real scheme of $U$ can be
described via $\langle \ \mathcal{C}_1|...|\mathcal{C}_k\rangle $,
where $|$ stands for a component homologous to
$(E_0)_{\mathbb{R}}$ and $\mathcal{C}_i$ encodes the arrangement
of the ovals in $Z_i$ (similar to the case of plane curves), for
each $i\in \{1,2,...,k\}$.

\begin{theorem} \label{sigma}{\rm(\cite{D2})}
A nonsingular real curve $U \in |3e_{\infty}|$ on $\Sigma_2$ is
determined up to rigid isotopy by its real scheme. There are 11
rigid isotopy classes, with real schemes $\langle \alpha |0\rangle
$, $1\leq \alpha \leq 4$, $\langle 0| \alpha \rangle $, $1\leq
\alpha \leq 4$, $\langle 0|0\rangle $, $\langle 1|1\rangle $, and
$\langle \mid \mid \mid \rangle $.
\end{theorem}

\begin{remark} \label{sigma'}
By analyzing the proof of Theorem~\ref{sigma}, one can easily see
that the curves with real schemes $\langle \alpha |0\rangle $ and
$\langle 0| \alpha \rangle $, $1\leq \alpha \leq 4$, are
isomorphic up to a real automorphism of $\Sigma_2$. Furthermore, a
stronger statement holds: any two pairs ($U$,$O$), where the real
scheme of $U$ is $\langle \alpha |0 \rangle $ with
$0\le\alpha\le3$ and $O$ is a distinguished oval of $U$, are
rigidly isotopic. For an alternative proof of Theorem~\ref{sigma}
and the last assertion, one can use the theory of the trigonal
curves, see~\cite{D3}.
\end{remark}

\subsection{Regular complete intersections of two real quadrics
in $\mathbb{P}_{\mathbb{R}}^4$}

The following is a special case of the rigid isotopy
classification of regular complete intersections of two quadrics
in $\mathbb{P}_{\mathbb{R}}^n$, due to S.~Lopez de
Medrano~\cite{LdM}.

\begin{theorem} \label{p4}{\rm(\cite{LdM})}
A regular complete intersection $Y$ of two real quadrics in
$\mathbb{P}_{\mathbb{R}}^4$ is determined up to rigid isotopy by
its real part $Y_{\mathbb{R}}$. There are seven rigid isotopy
classes, with $Y_{\mathbb{R}}=V_6$, $V_4$, $V_2$, $S_1$, $2S$,
$S$, or $\varnothing$.
\end{theorem}

\subsection{Real root schemes} Let $Z=\Sigma_k$, $k\geq 0$, with
the standard real structure. Since we use $\Sigma_2$ and
$\Sigma_4$ in this paper we will consider only the cases $k=2n$.
For $k=2n+1$ and further details, see~\cite{D2}. Consider a real
curve $U\in |2e_\infty+pl|$, $p\geq 0$, and a real curve
$Q=E_0\cup F$, where $E_0$ is the exceptional section and $F\in
|e_\infty|$ is a generic real section of $Z$. The complement
$Z_{\mathbb{R}}\setminus Q_{\mathbb{R}}$ consists of two connected
orientable components. Fix one of them and let $Z^-$ denote its
closure. Fix an orientation of $F_{\mathbb{R}}\subset \partial
Z^-$. Assume that $U$ does not contain any generatrix of $Z$, is
transversal to $F$ and $U_{\mathbb{R}}$ lies entirely in $Z^-$.
Fix an auxiliary real generatrix $L$ of $Z$ transversal to $U\cup
E_0$. Consider a real coordinate system $(x,y)$ in the affine part
$Z\setminus (E_0\cup L)$ whose $x$-axis is $F$. Choose the
positive direction of the $y$-axis so that the upper half-plane
lies in $Z^-$. In these coordinates $U$ has equation
$a(x)y^2+b(x)y+c(x)=0$, where $a,b$, and $c$ are real polynomials
of degree $p$, $p+k$, and $p+2k$, respectively. Let
$\Delta=b^2-4ac$ and let $\mu(x)$ and $\nu(x)$ denote the
multiplicity of a point $x\in F$ in $a$ and $\Delta$,
respectively. Consider the sets $$\mathcal{A}_{\mathbb{R}}=\{x\in
F_{\mathbb{R}}\ |\ \mu(x)\geq 1\},\ \ \ \mathcal{A}=\{x\in F\ | \
\mu(x)\geq 1\},$$
$$\mathcal{D}_{\mathbb{R}}=\{x\in F_{\mathbb{R}}\ |\ \Delta(x)\geq 0\},
\ \ \ \mathcal{D}_r=\{x\in F\ |\ \nu(x)\geq r\},\ r\geq 1, \ \ \
\mathcal{D}=\mathcal{D}_2\cup \mathcal{D}_{\mathbb{R}}.$$ The
multiplicity functions $\mu$ and $\nu$ are invariant under complex
conjugation. Identify $F$ with the base $B\cong \mathbb{P}^1$ of
the ruling of $Z$. Thus, $B_{\mathbb{R}}$ receives an orientation,
$\mathcal{A}$ and $\mathcal{D}$ can be regarded as subsets of $B$,
and, $\mu$ and $\nu$ are functions defined on $B$. The
\textit{root marking} of $(U,Q)$ is the triple
$(\mathcal{B},\mathcal{D},\mathcal{A})$ equipped with the complex
conjugation in $B$ and the following structures:
\begin{enumerate}
\item the orientation of $B_{\mathbb{R}}$; \item the multiplicity
functions $\mu$ and $\nu$.
\end{enumerate}
An \textit{isotopy} of root markings is an equivariant isotopy of
triples $(\mathcal{B},\mathcal{D},\mathcal{A})$ followed by a
continuous change of the orientation of $B_{\mathbb{R}}$, $\mu$,
and $\nu$ restricted to $\mathcal{D}$. A \textit{root scheme} is
an equivalence class of root markings up to isotopy. The
\textit{real root marking} of $(U,Q)$ is the triple
$(\mathcal{B}_{\mathbb{R}},\mathcal{D}_{\mathbb{R}},\mathcal{A}_{\mathbb{R}})$
equipped with (1) and (2) above. A \textit{real root scheme} is an
equivalence class of real root markings up to isotopy.

\begin{theorem} {\rm(\cite{D2})}\label{Sig4}
Let $Z=\Sigma_{4}$ (with the standard real structure), let $U\in
|2e_{\infty}|$ be a nonsingular real curve on $Z$, let $F\in
|e_{\infty}|$ be a generic real section transversal to $U$, and
let $E_0$ be the exceptional section. If $U_{\mathbb{R}}$ belongs
to the closure of one of the two components of
$Z_{\mathbb{R}}\setminus ((E_0)_{\mathbb{R}}\cup F_{\mathbb{R}})$,
then, up to rigid isotopy and automorphism of $Z$, the pair
$(U,F)$ is determined by its real root scheme or, equivalently, by
the real scheme of $U$. The latter consists either of $a=0,...,4$
ovals ({\em i.e.}, components bounding disks) or of two components
isotopic to $F_{\mathbb{R}}$.
\end{theorem}

\subsection{Suitable pairs.}\label{suip}
Let $U\in |2e_{\infty}+2l|$ be a
reduced (does not contain any multiple component) real curve on
$\Sigma_{2}$ with the standard real structure. Assume that $U$ is
nonsingular outside of $E_0$ and does not contain $E_0$ as a
component. Then $U$ and $E_0$ intersect with multiplicity $2$. So $U$
either intersects $E_0$ transversally at two points, or is tangent
to it at one point, or has a single singular point of type
$\textbf{A}_{r-2},\ r\geq 3$, on $E_0$; the \textit{grade} of $U$
is said to be $1$, $2$, or $r$ respectively. A curve $U$ as above
is called \textit{suitable} if either its grade is even or grade
is odd and the two branches of $U$ at $E_0$ are conjugate to each
other. A pair $(U,F)$ is called a \textit{suitable pair} if $U$ is a
suitable curve and $F\in |e_{\infty}|$ a nonsingular real section
transversal to $U$ such that $U_{\mathbb{R}}$ belongs to the
closure of a single connected component of
${(\Sigma_2)}_{\mathbb{R}}\setminus ((E_0)_{\mathbb{R}}\cup
F_{\mathbb{R}})$. The \textit{grade} of a suitable pair $(U,F)$ is
the grade of $U$. The condition that $U_{\mathbb{R}}$ should
belong to the closure of a single connected component of
${(\Sigma_2)}_{\mathbb{R}}\setminus ((E_0)_{\mathbb{R}}\cup
F_{\mathbb{R}})$ guarantees that the real \textit{DPN}-double
$(Y,B)$ ({\em i.e.}, the resolution of singularities of the double
covering of $\Sigma_2$ branched over $U$, where the rational
components of $B$ correspond to $E_0$ and the irrational component
of $B$ corresponds to $F$) of $(\Sigma_{2};U,E_0\cup F)$, where
$(U,F)$ is a suitable pair, corresponds to a real Enriques surface
by inverse Donaldson's trick.

All the pairs $(U,F)$ satisfying the hypothesis of the following
theorem are suitable.

\begin{theorem} {\rm(\cite{D2})}\label{SP}
Let $Z=\Sigma_{2}$ (with the standard real structure), let $U\in
|2e_{\infty}+2l|$ be a reduced real curve on $Z$, nonsingular
outside the exceptional section $E_0$ and not containing $E_0$ as
a component, and let $F\in |e_{\infty}|$ be a generic real section
transversal to $U$. If $U_{\mathbb{R}}$ belongs to the closure of
a single connected component of $Z_{\mathbb{R}}\setminus
((E_0)_{\mathbb{R}}\cup F_{\mathbb{R}})$, then, up to rigid
isotopy and automorphism of $Z$, the pair $(U,F)$ is determined by
its real root scheme or, equivalently, by the type of the singular
point of $U$ (if any) and the topology of the pair
$(Z_{\mathbb{R}},U_{\mathbb{R}}\cup (E_0)_{\mathbb{R}})$.
\end{theorem}

\begin{table} [h]
\caption{Real root schemes of some curves $U\in |2e_{\infty}+pl|$
on $\Sigma_{2k}$} \catcode`\@\active

\def\.{@(\bullet)}
\def\*{@(\circ)}
\def\-{@---}
\def\={@\---}
\def\w#1{\*@*<|><\tt#1|>}
\def\b#1{\.@*<|><\tt#1|>}
\def\t#1{\mathord{\scriptstyle\times}#1}
\def\bt#1{\overline{\t{#1}}}
\def\tt#1{\t{\vcenter{\hbox{$\scriptstyle#1$}}}}
\let\0\varnothing
\def\(#1)\[#2]\[#3]{@\l[\quad$\{#2\}\+\{#3\}$]}
\def\[#1]{@\l[\let\b\BAR$#1$\qquad]}
\def\[#1]{}
\def\c{0}
\def\n{2}
\def\>{@[$\!\longrightarrow\!\!\!$]}
\def\\{\strut\CR}
\hbox to\hsize{\hss \DDbox \expand20@\l[\quad Real root scheme]
  \($\CP$)\[\er1]\[\er2]\\
\noalign{\smallskip\hrule\smallskip}
\[\0]
  \-\.\=\.\-\.\=\.\-\.\=\.\-\.\=\.\-&&
  \()\[(V_4\+S)\+(\0)]\[(2S)\+(2S)]\\
\[(\bA_0,\b\bA_0)]
  \-\.\=\.\-\.\=\.\-\.\=\.\-\.\=\.\-&&
  \()\[(V_3\+V_1)\+(\0)]\[(2S)\+(2S)]\\
\[(\bA_1)]
  \-\.\=\*\=\.\-\.\=\.\-\.\=\.\-\.\=\.\-
  \()\[(V_3\+S)\+(\0)]\[(V_1\+S)\+(2S)]\\
\[(\bB_2)]
  \-\w2\-\.\=\.\-\.\=\.\-\.\=\.\-&&&&
  \()\[(V_3\+S)\+(V_1)]\[(2S)\+(S)]\\
\[(\bB_3)]
  \-\w3\=\.\-\.\=\.\-\.\=\.\-&&&&&&
  \()\[(V_3\+S)\+(S)]\[(V_1\+S)\+(S)]\\
\[(\bB_4)]
  \-\w4\-\.\=\.\-\.\=\.\-&&&&&&&&
  \()\[(V_3\+V_1\+S)\+(S)]\[(S)\+(S)]\\
\endDD\hss}

\eightpoint
\bigskip
\par
\leftskip0pt \rightskip0pt \textbf{Comments:} The first column
indicates the real root schemes of pairs $(U,F)$ and the second
column indicates the quarter decomposition of the real part
$E_{\mathbb{R}}$ of the real Enriques surfaces obtained from
$(\Sigma_{2k};U,E_0\cup F)$. For the first row $p=0$ and $k=2$ and
for the others $p=2$ and $k=1$. In the schemes,
$\DDcenter\.\endDD\ $ represents a real root of~$\Delta$ and
$\DDcenter\*\endDD\ $ represents a real root of~$a$ (necessary
2-fold), that corresponds to the real intersection point of $U$
and $E_0$. The number over a $\DDcenter\*\endDD\ $-vertex
indicates the multiplicity of the corresponding root in~$\Delta$
(when greater than~$1$). The segments $\ \DDcenter\.\=\.\endDD\ $
correspond to ovals of $U_{\mathbb{R}}$. Only extremal root
schemes are listed; the others are obtained by removing one or
several segments $\ \DDcenter\.\=\.\endDD\ $.

\label{root}
\end{table}

In Table~\ref{root}, we list the extremal real root schemes of
some pairs $(U,F)$ mentioned in Theorem~\ref{Sig4} and
Theorem~\ref{SP} that are used in the proof of the main result.
The complete lists can be found in~\cite{D2}.

\begin{remark}\label{common}
Each real root marking gives rise to a connected family of pairs
$(U,Q)$ such that there is a bijection between the ovals of each
curve $U$ and the segments of the real root marking. Recall that
these curves are defined by explicit equations. Then both
Theorems~\ref{Sig4} and~\ref{SP} can be refined as follows:
\begin{enumerate}
\item Each isotopy of real root markings is followed by a rigid
isotopy of curves that is consistent with the bijection between
ovals and segments. \item Any symmetry of a real root marking (not
necessarily preserving the orientation of $B_{\mathbb{R}}$) is
induced by an automorphism of $\Sigma_{2k}$, $k\ge0$, preserving
appropriate pairs $(U,F)$ and consistent with the bijection
between ovals and segments.
\end{enumerate}
\end{remark}

\section{Reduction to \textit{DPN}-pairs}\label{DPN}

\subsection{Donaldson's trick}\label{Don}
The equivariant version of Donaldson's trick employs the
hyper-K\"ahler structure to change the complex structure of the
covering $K3$-surface $X$ so that $t^{(1)}$ is holomorphic, and
$t^{(2)}$ and $\tau$ are anti-holomorphic. Furthermore,
$Y=X/t^{(1)}$ is a real rational surface, where the real structure
is the common descent of $\tau$ and $t^{(2)}$, and $B\cong \Fix
t^{(1)}$ is a nonsingular curve on $Y$. As a result, the problem
about real Enriques surfaces is reduced to the study of certain
auxiliary objects, like real plane quartics, space cubics,
intersections of two quadrics in $\mathbb{P}_{\mathbb{R}}^4$, etc.

\begin{theorem}\label{Donaldson} {\rm(\cite{D1})}
Donaldson's construction establishes a one-to-one correspondence
between the set of deformation classes of real Enriques surfaces
with distinguished nonempty half ({\em i.e.}, pairs $(E,
E_{\mathbb{R}}^{(1)})$ with $E_{\mathbb{R}}^{(1)} \neq
\varnothing$) and the set of deformation classes of pairs $(Y,B)$,
where $Y$ is a real rational surface and $B \subset Y$ is a
nonsingular real curve such that
\begin{enumerate}
\item $B$ is anti-bicanonical, \item the real point set of $B$ is
empty, and \item $B$ is not linked with the real point set
$Y_{\mathbb{R}}$ of $Y$.
\end{enumerate}
One has $E_{\mathbb{R}}^{(2)}=Y_{\mathbb{R}}$ and
$E_{\mathbb{R}}^{(1)}=B/t^{(2)}$.
\end{theorem}

(A real curve $B \subset Y$ with $B_{\mathbb{R}}= \varnothing$ is
said to be {\em not linked with $Y_{\mathbb{R}}$} if for any path
$\gamma : [0,1] \rightarrow Y\setminus B$ with $\gamma (0), \gamma
(1) \in Y_{\mathbb{R}}$, the loop $\gamma^{-1} \cdot \conj_Y
\gamma$ is $\mathbb{Z}/2$-homologous to zero in $Y\setminus B$.)

In the above theorem, the first condition on $B$ guarantees that
the double covering $X$ of $Y$ branched over $B$ is a
$K3$-surface; and the other two conditions ensure the existence of
a fixed point free lift of the real structure on $Y$ to $X$,
see~\cite{D1}. The statement deals with deformation classes rather
than individual surfaces because the construction involves a
certain choice (that of an invariant K\"ahler class).

\subsection{Inverse Donaldson's trick} \label{inverse}
Since we want to construct deformation families of real Enriques
surfaces with particular properties, we are using Donaldson's
construction backwards. Strictly speaking, Donaldson's trick is
not invertible. However, it establishes a bijection between the
sets of deformation classes (see Theorem \ref{Donaldson}); thus,
at the level of deformation classes one can speak about 'inverse
Donaldson's trick'.

Before explaining the construction, recall some properties of
$K3$-surfaces. Let $a$ be a holomorphic involution of a
$K3$-surface $X$ equipped with the complex structure defined by a
holomorphic form $\omega $. Then there are three possibilities for
the fixed point set $\Fix a$ of $a$:
\begin{enumerate}
\item it may be empty, or \item It may consist of isolated points,
or \item it may consist of curves.
\end{enumerate}
The following is straightforward:
\begin{enumerate}
\item if $\dim_{\mathbb{C}}\Fix a=0$, then $a^*\omega=\omega$,
\item if $\dim_{\mathbb{C}}\Fix a=\pm 1$, then
$a^*\omega=-\omega$.
\end{enumerate}

Let $\conj $ be a real structure on $X$. Then ${\conj}
^*\omega=\lambda \overline{\omega}$ for some $\lambda \in
\mathbb{C}^*$. Clearly, $w$ can be chosen (uniquely up to real
factor) so that ${\conj}^*\omega=\overline{\omega}$. We always
assume this choice and we denote by $\Re{\omega}$ and
$\Im{\omega}$ the real part $(\omega+\overline \omega)/2 $ and the
imaginary part $(\omega-\overline \omega)/2$ of $\omega$,
respectively.

Let $Y$ be a real rational surface with a nonsingular
anti-bicanonical real curve $B \subset Y$ such that
$B_{\mathbb{R}}=\varnothing $ and $B$ is not linked with the real
point set $Y_{\mathbb{R}}$ of $Y$. Let $X$ be the (real) double
covering $K3$-surface branched over $B$, $\widetilde{p}:X
\rightarrow Y$ the covering projection and $\phi:X \rightarrow X$
the deck translation of $\widetilde{p}$. Then $\phi$ is a
holomorphic involution with nonempty fixed point set. There exist
two liftings $c^{(1)},c^{(2)}: X \rightarrow X$ of the real
structure $\conj: Y \rightarrow Y$ to $X$, which are both
anti-holomorphic involutions. They commute with each other and
with $\phi$, and their composition is $\phi$. Because of the
requirements on $B$, at least one of these involutions is fixed
point free. Assume that it is $c^{(1)}$.

Pick a holomorphic 2-form $\mu$ with the real and imaginary parts
$\Re{\mu}, \Im{\mu}$, respectively, and a fundamental K\"ahler
form $\nu$. Due to the Calabi-Yau theorem, there exists a unique
K\"ahler-Einstein metric with fundamental class $[\nu]$,
see~\cite{H}. After normalizing $\mu$ so that
$\Re{\mu}^2=\Im{\mu}^2={\nu}^2=2\Vol X$, we get three complex
structures on $X$ given by the forms: $$ \mu=\Re{\mu}+i\Im{\mu},\
\ \ \widetilde{\mu}=\nu+i\Re{\mu} ,\ ~\textmd{and}~\ \Im{\mu}+i\nu
.$$ Let $\widetilde{X}$ be the surface $X$ equipped with the
complex structure defined by $\widetilde{\mu}$. Since $c^{(1)}$ is
an anti-holomorphic involution of $X$, the holomorphic form $\mu$
and the fundamental K\"ahler form $\nu$ can be chosen so that
$(c^{(1)})^*\mu =-\overline{\mu}$ and $(c^{(1)})^*\nu=-\nu$. Then
$(c^{(1)})^*\widetilde{\mu}=-\widetilde{\mu}$ and, hence,
$c^{(1)}$ is holomorphic on $\widetilde{X}$. Since $\phi$ is a
holomorphic involution of $X$ commuting with $c^{(1)}$, $\
\phi^*\mu=-\mu$ and $\nu$ can be chosen $\phi^*$-invariant so that
$\phi^*\widetilde{\mu}=\overline{\widetilde{\mu}}$, {\em i.e.},
the involution $\phi$ is anti-holomorphic on $\widetilde{X}$. Then
$E=\widetilde{X}/c^{(1)}$ is a real Enriques surface (the real
structure being common descent of $\phi$ and $c^{(2)}$) and the
projection $p:\widetilde{X} \rightarrow E$ is a real double
covering.  Hence we have $Y_{\mathbb{R}}=E_{\mathbb{R}}^{(2)}$ and
$B/c^{(2)}=E_{\mathbb{R}}^{(1)}$.

\subsection{The case of Del Pezzo surfaces} The deformation classification of real Enriques
surfaces with a distinguished half $E_\mathbb{R}^{(1)}=V_{d+2}$,
$d\geq 1$ is reduced to that of real unnodal Del Pezzo surfaces of
degree $d$, $d\geq 1$, with a nonsingular anti-bicanonical curve
$B\cong S_{g}$, $g\geq 2$.

\begin{theorem} {\rm(\cite{D2})}\label{E-DP}
There is a natural surjective map from the set of deformation
classes of real unnodal Del Pezzo surfaces Y of degree $d$, $d\geq
1$, onto the set of deformation classes of real Enriques surfaces
with $E_\mathbb{R}^{(1)}=V_{d+2}$, $d\geq 1$. Under this
correspondence $Y_\mathbb{R}=E_\mathbb{R}^{(2)}$ and
$Y/\conj=E/\conj$.
\end{theorem}

\begin{remark}
In fact, the correspondence is bijective.
\end{remark}

Proof of Theorem~\ref{E-DP} reduces, mainly, to showing that a
generic deformation of unnodal Del Pezzo surfaces~$Y_t$ can be
extended to a deformation of pairs $(Y_t,B_t)$, where $B_t\subset
Y_t$ are real anti-bicanonical curves satisfying the hypotheses of
Theorem~\ref{Donaldson}. This gives a deformation of the covering
$K3$-surfaces. Then it remains to choose a continuous family of
invariant K\"ahler metrics, and inverse Donaldson's trick applies.
Thus, the following stronger result holds.

\begin{theorem} \label{delpezzo}
A generic deformation of real unnodal Del Pezzo surfaces~$Y$ of
degree $d$, $d\geq1$, defines a deformation of real Enriques
surfaces with $E_\mathbb{R}^{(1)}=V_{d+2}$, $d\geq1$, obtained
from~$Y$ by inverse Donaldson's trick.
\end{theorem}

\subsection{The case of $(2,r)$ surfaces} The deformation classification of real Enriques surfaces with
disconnected $E_{\mathbb{R}}^{(1)}=V_{3}\sqcup ...$ is reduced to
that of real $(2,r)$-surfaces, $r\geq 1$ with a real nonsingular
anti-bicanonical curve $B\cong S_{2}\sqcup rS$ and, hence, to the
rigid isotopy classification of suitable pairs.

\begin{lemma} {\rm(\cite{D2})} \label{supa}
There is a natural surjective map from the set of rigid isotopy
classes of suitable pairs of grade $r$ onto the set of deformation
classes of real Enriques surfaces with
$E_{\mathbb{R}}^{(1)}=V_{3}\sqcup \frac{r}{2}S$, if $r$ is even,
or $E_{\mathbb{R}}^{(1)}=V_{3}\sqcup V_{1} \sqcup \frac{r-1}{2}S$,
if $r$ is odd.
\end{lemma}

Proof of the above lemma is based on showing that a generic rigid
isotopy of suitable pairs $(U_s,F_s)$ defines a deformation of the
\textit{DPN}-doubles $(Y_s,B_s)$ of $(\Sigma_2;U_s,E_0\cup F_s)$,
so a deformation of the covering $K3$-surfaces. Then it remains to
choose a continuous family of invariant K\"ahler metrics, to
obtain a deformation of the corresponding real Enriques surfaces
obtained by inverse Donaldson's trick which implies the following
stronger result.

\begin{theorem}\label{gensp}
A generic rigid isotopy of suitable pairs $(U,F)$ of grade $r$
defines a deformation of the real Enriques surfaces with
$E_{\mathbb{R}}^{(1)}=V_{3}\sqcup \frac{r}{2}S$, if $r$ is even,
or $E_{\mathbb{R}}^{(1)}=V_{3}\sqcup V_{1} \sqcup \frac{r-1}{2}S$,
if $r$ is odd.
\end{theorem}

\section{Lifting Involutions}\label{Lifting}
Let $Z$ be a simply connected surface and $\pi :Y\rightarrow Z$ a
branched double covering with the branch divisor $C$. Then any
involution $a:Z\rightarrow Z$ preserving $C$ as a divisor admits
two lifts to $Y$, which commute with each other and with the deck
translation of the covering. If $\Fix a\neq \varnothing $, then
both lifts are also involutions. Any fixed point of $a$ in
$Z\setminus C$ has two pullbacks on $Y$. One of the lifts fixes
these two points and the other one permutes them.

In this section we will use the notation of subsection
\ref{inverse}.

\begin{lemma} \label{lift1}
Let $Z$ be a real quadric cone in $\mathbb{P}^3$, let $C \subset
Z$ be a nonsingular real cubic section disjoint from the vertex,
and let $a:Z\to Z$ be an involution preserving~$C$ and such that
$\Fix{a}\cap C \neq \varnothing$. Then $a$ lifts to four distinct
involutions on the covering $K3$-surface~$X$ and at least one of
the four lifts defines an automorphism of an appropriate real
Enriques surface obtained from~$X$ by inverse Donaldson's trick.
\end{lemma}

\begin{proof}
According to the models of Del Pezzo surfaces \cite{Dem}, the
double covering of $Z$ branched at the vertex and over $C$ is a
real unnodal Del Pezzo surface $Y$ of degree $d=1$. The pullback
$\widetilde{p}\in Y$ of any point $p \in \Fix{a} \cap C$ is a
fixed point of any lift of~$a$ to~$Y$. Let $p' \in
\Fix{a}\setminus C$ be in a small neighborhood of $p$. Then $p'$
has two pullbacks $p_1$ and $p_2$ in $Y$. Let $a_1$ be the lift of
$a$ to $Y$ that permutes $p_1$ and $p_2$. Then $\widetilde{p}$ is
an isolated fixed point of~$a_1$. (Note that we do \textbf{not}
assert that \textbf{all} fixed points of~$a_1$ are isolated.) Pick
an $a_1$-invariant admissible branch curve $B\subset Y$ with
$\widetilde{p} \notin B$. Denote by~$X$ the double covering of~$Y$
branched over~$B$ and by $a_2$, the lift of $a_1$ to $X$ that
fixes the two pullbacks of $\widetilde{p}$. Then the pullbacks of
$\widetilde{p}$ are isolated fixed points of $a_2$. Since $X$ is a
$K3$-surface, $\Fix{a_2}$ consists of isolated points only, and
$(a_2)^*\mu=\mu$. We can choose for $\nu$ a generic fundamental
K\"ahler form preserved by $\phi$, $c^{(1)}$, $c^{(2)}$, and
$a_2$. Then we have $(a_2)^*\widetilde{\mu}=\widetilde{\mu}$, {\em
i.e.}, $a_2$ is also holomorphic on $\widetilde{X}$. With the
projection $p:\widetilde{X} \rightarrow E$, $a_2$ defines an
automorphism $\widetilde{a}$ of $E$.
\end{proof}

\begin{lemma} \label{lift2}
Let $Y$ be a real unnodal Del Pezzo surface of degree $d=2$ with
$Y_{\mathbb{R}}=2V_1$ and let $\Gamma$ be the deck translation
involution of the double covering $Y\rightarrow \mathbb{P}^2$
whose branch locus is a nonsingular quartic $C$ with
$C_{\mathbb{R}}=\varnothing$. Then $\Gamma$ lifts to two distinct
involutions on the covering $K3$-surface~$X$, and one of the lifts
defines an automorphism of an appropriate real Enriques surface
obtained from~$X$ by inverse Donaldson's trick.
\end{lemma}

\begin{proof}
Pick a $\Gamma$-invariant admissible branch curve $B\subset Y$ and
denote by~$X$ the double covering $K3$-surface of~$Y$ branched
over~$B$. Due to the adjunction formula, $\Fix {\Gamma} \cap B
\neq \varnothing$. Hence, as in the previous case, we can choose a
lift $a$ of $\Gamma$ to $X$ having isolated fixed points. Since
$X$ is a $K3$-surface, $\Fix{a}$ consists of isolated points only,
and $a^*\mu=\mu$. We can choose for $\nu$ a generic fundamental
form preserved by $\phi$, $c^{(1)}$, $c^{(2)}$, and $a$. Then we
have $a^*\widetilde{\mu}=\widetilde{\mu}$, {\em i.e.}, $a$ is also
holomorphic on $\widetilde{X}$. With the projection
$p:\widetilde{X} \rightarrow E$, $a$ defines an automorphism
$\widetilde{a}$ of $E$.
\end{proof}

\begin{lemma} \label{lift3}
Let $Z=\Sigma_4$ (with the standard real structure), and $U\in
|2e_{\infty}|$ a nonsingular real curve. Let $a:Z\rightarrow Z$ be
an involution preserving $U$ and such that $\Fix a\cap U\neq
\varnothing$. Then $a$ lifts to four distinct involutions on the
covering $K3$-surface~$X$ and at least one of the four lifts
defines an automorphism of an appropriate real Enriques surface
obtained from~$X$ by inverse Donaldson's trick.
\end{lemma}

\begin{proof}
For a nonsingular real curve $F\in |e_{\infty}|$ in $Z$, if
$U_{\mathbb{R}}$ is contained in a connected component of
$Z_{\mathbb{R}}\backslash((E_0)_\mathbb{R}\sqcup F_{\mathbb{R}})$
then the \textit{DPN}-double $(Y,B)$ of $(Z;U,E_0\cup F)$ is as
follows: $Y$ is a real unnodal $(3,2)$-surface, and, $B$ is an
admissible branch curve with two rational components which are
conjugate to each other and $[B]=0$ in $H_2(X)$ where $X$ is the
covering $K3$-surface of $Y$ branched over $B$ (see \cite{D2}).
Any point $p\in \Fix{a}\cap U$ has a unique pullback
$\widetilde{p}\in Y$ which is a fixed point of both lifts of $a$
to $Y$. Any point $p'\in \Fix{a}\setminus U$, in a small
neighborhood of $p$, has two pullbacks $p_1$ and $p_2$ in $Y$. If
$a_1$ is the lift of $a$ to $Y$ that permutes $p_1$ and $p_2$ then
$\widetilde{p}$ is an isolated fixed point of $a_1$. Choose $F\in
|e_{\infty}|$ and the point $p\in \Fix{a}\cap U$ in such a way
that $B$ is $a_1$-invariant and $\widetilde{p}\notin B$. Let $X$
be the double covering of~$Y$ branched over~$B$ and let $a_2$ be
the lift of $a_1$ to $X$ that fixes the two pullbacks of
$\widetilde{p}$. Then the pullbacks of $\widetilde{p}$ are
isolated fixed points of $a_2$. Since $X$ is a $K3$-surface,
$\Fix{a_2}$ consists of isolated points only, and
$(a_2)^*\mu=\mu$. The result follows by making the same choices as
in the proof of Lemma~\ref{lift1}.
\end{proof}

\begin{lemma} \label{lift4}
Let $Z=\Sigma_2$ (with the standard real structure), let $U\in
|2e_\infty+2l|$ be a suitable curve on $Z$, and let
$a:Z\rightarrow Z$ be an involution preserving $U$ such that $\Fix
a\cap U\neq \varnothing$. Then $a$ lifts to four distinct
involutions on the covering $K3$-surface~$X$ and at least one of
the four lifts defines an automorphism of an appropriate real
Enriques surface obtained from~$X$ by inverse Donaldson's trick.
\end{lemma}

\begin{proof}
For a nonsingular real section $F\in |e_{\infty}|$ in $Z$, if
$(U,F)$ is a suitable pair then the \textit{DPN}-double of
$(Z;U,E_0\cup F)$ is $(Y,B)$ where $Y$ is a $(2,r)$-surface and
$B$ is an admissible branch curve on $Y$ (see~\cite{D2}). Thus,
for any such curve $F$, we can make choices of the points $p\in
\Fix{a}\cap U$ and $p'\in \Fix{a}\setminus U$, and the lift $a_1$
of $a$ to $Y$ in the same way that we did in the proof of
Lemma~\ref{lift3} so that $\widetilde{p}\in Y$ will be an isolated
fixed point of $a_1$. Choose $F\in |e_{\infty}|$ and the point
$p\in \Fix{a}\cap U$ in such a way that $B$ is $a_1$-invariant and
does not contain $\widetilde{p}$. Then the result follows by
making the same choices as in the proof of Lemma~\ref{lift3}.
\end{proof}

\section{Main results}\label{Main}

\begin{theorem} \label{DP}
With one exception, any permutation of homeomorphic components of
the half $E_{\mathbb{R}}^{(2)}$ of a real Enriques surface with a
distinguished half $E_{\mathbb{R}}^{(1)}= V_{d+2}$, $d\geq 1$, can
be realized by deformations and automorphisms. In the exceptional
case $E_{\mathbb{R}}= \{V_3 \}\sqcup \{V_1 \sqcup 4S \}$, the
realized group is $\mathbb{Z}_2 \times \mathbb{Z}_2\subset S_4$.
\end{theorem}

\begin{remark}\label{no.others}
In the exceptional case, the Pontrjagin-Viro form (see~\cite{D4})
is well defined. It defines a decomposition of $E_{\mathbb{R}}$
into quarters, which is a topological invariant. The decomposition
of~$E_\mathbb{R}^{(2)}$ is $(V_1\sqcup 2S)\sqcup (2S)$. Obviously,
one cannot permute the spheres belonging to different quarters
(even topologically), and Theorem~\ref{DP} states that a
permutation of the spherical components can be realized if and
only if it preserves the quarter decomposition.
\end{remark}

\begin{proof}
The deformation classification of real Enriques surfaces with a
distinguished half $E_\mathbb{R}^{(1)}=V_{d+2}$, $d\geq1$ is
reduced to that of real unnodal Del Pezzo surfaces of degree $d$,
$d\geq1$, with a nonsingular anti-bicanonical curve $B\cong
S_{g}$, $g\geq2$ (see~\cite{Dem} for the models of Del Pezzo
surfaces). It always suffices to construct a particular surface
(within each deformation class) that has a desired automorphism or
`auto-deformation'. We proceed case by case. Among the extremal
types listed in~\cite{D2}, we need to consider only the following
types (as in the other cases there are no homeomorphic components)
and all their derivatives $(E_{\mathbb{R}}^{(1)},\cdot )$ obtained
from the extremal ones by sequences of topological Morse
simplifications of $E_\mathbb{R}^{(2)}$:

\begin{enumerate}
\item $E_{\mathbb{R}}^{(1)}=V_3$;\quad $E_{\mathbb{R}}^{(2)}=V_1
\sqcup 4S$; \item $E_{\mathbb{R}}^{(1)}=V_4$;\quad
$E_{\mathbb{R}}^{(2)}=2V_1$; \item
$E_{\mathbb{R}}^{(1)}=V_4$;\quad $E_{\mathbb{R}}^{(2)}=4S$; \item
$E_{\mathbb{R}}^{(1)}=V_6$;\quad $E_{\mathbb{R}}^{(2)}=2S$.
\end{enumerate}
\vspace{0.3cm}

\textbf{Case 1}: Here we consider the 3 subcases:
$$E_{\mathbb{R}}=\{V_3 \} \sqcup \{V_1 \sqcup iS\},~~i=2,3,4.$$
The corresponding surface $Y$ obtained by Donaldson's trick (see
subsection \ref{Donaldson}) is a real unnodal Del Pezzo surface of
degree 1 with $Y_{\mathbb{R}}=V_1 \sqcup iS$, $i=2,3,4$. According
to the models of Del Pezzo surfaces, the anti-bicanonical system
$|{-2K}|$ maps $Y$ onto an irreducible singular quadric (cone) $Z$
in $\mathbb{P}^3$. This map $\varphi :Y \rightarrow Z$ is of
degree 2 and its branch locus consists of the vertex $V$ of $Z$
and a nonsingular cubic section $C$ disjoint from $V$ whose real
part $C_{\mathbb{R}}$ consists of $i$~ovals and a component
noncontractible in $Z_{\mathbb{R}}\setminus\{V\}$.
The real part $Y_{\mathbb{R}}$ is the double covering of the domain $D$
consisting of $i$~disks bounded by the ovals of $C_{\mathbb{R}}$
and of the part of $Z_{\mathbb{R}}$ bounded by the noncontractible
component of $C_{\mathbb{R}}$ and $V$. The map $\varphi$ lifts to
a degree $2$ map $\tilde{\varphi} :Y \rightarrow
Z=\Sigma_{2}\subset \mathbb{P}^3\times \mathbb{P}^1$ ($\Sigma_2$
with standard real structure). The branch set of $\tilde{\varphi}$
is the union of $E_0$ and a real nonsingular curve $C' \in
|3e_{\infty}|$. Rigid isotopy class of $C$ is induced by that of
$C'$. From Theorem~\ref{sigma} and Remark~\ref{sigma'}, for
each~$i=2,3$ and $4$, there is one rigid isotopy class of $C$ up
to isomorphism.

Clearly, a rigid isotopy of~$C$ in~$Z$ defines a deformation
of~$Y$, and an auto-involution of~$Z$, preserving~$C$ and having
nonempty fixed point set, lifts to an involution on~$Y$. Thus, in
view of Theorem~\ref{delpezzo} and Lemma~\ref{lift1}, it suffices
to realize certain permutations of the ovals of a particular curve
(in each rigid isotopy class)~$C$ by rigid isotopies and/or
involutive automorphisms of~$Z$ (in the latter case taking care
that the fixed point set of the involution intersects~$C$).

\begin{figure}[h]
\centering \epsfig{file=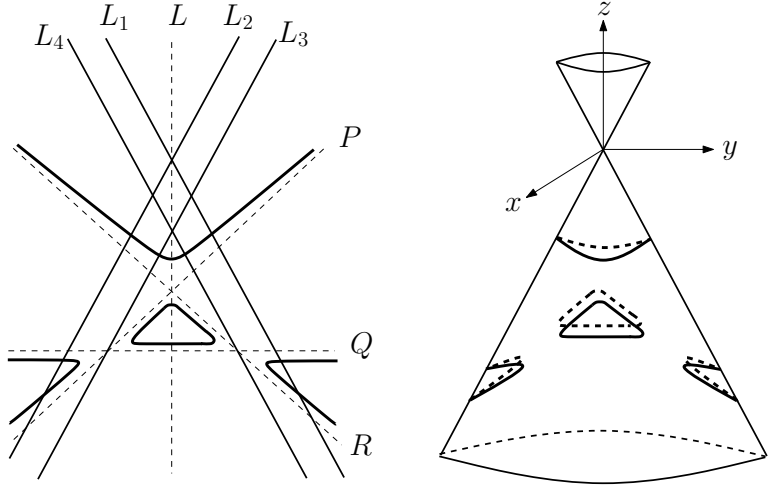,width=4in} \caption{Elements of the
construction of a quadric cone $Z\subset \mathbb{P}^3$ and a
symmetric cubic section $C\subset Z$ (left), and an example of the
maximal case (right)} \label{cubic}
\end{figure}

For each~$i=2,3$ and $4$, let $C=Z\cap S$, where $Z$ and $S$ in
$\mathbb{P}^3$ are constructed (due to S. Finashin,
see~\cite{Fin}) as follows: Let $Z$ be the quadric cone that is
the double covering of the plane branched over $L_3$ and $L_4$ if
$i=2$, $L_1$ and $L_3$ if $i=3$, and, $L_1$ and $L_2$ if $i=4$
(see Figure~\ref{cubic}). Let $S$ be the cubic surface that is the
pullback of the cubic curve, which is symmetric with respect to
the line $L$, and is obtained by a perturbation of the lines $P$,
$Q$ and $R$ (dotted lines, see Figure~\ref{cubic}). For $i=2$, the
symmetry of the cone with respect to the $yz$-plane permutes the
ovals of $C$. For $i=3$, it suffices to permute one pair of ovals,
see Remark~\ref{sigma'}, and the symmetry of the cone with respect
to the $yz$-plane does permute the opposite ovals of $C$. For
$i=4$, the symmetries of the cone with respect to the $yz$-plane
and $xz$-plane permutes the opposite ovals of $C$. Fixed point set
of each symmetry intersects $C$. Thus, we obtain the groups $S_2$,
$S_3$ and $\mathbb{Z}_2 \times \mathbb{Z}_2\subset S_4$ for
$i=2,3$ and $4$, respectively. For $i=4$, the fact that other
permutations cannot be realized is explained in
Remark~\ref{no.others}.

In cases 2 and 3 below, the corresponding surface $Y$ obtained by
Donaldson's trick is a real unnodal Del Pezzo surface of
degree~$2$. According to the models of Del Pezzo surfaces, the
anti-canonical system of $Y$ defines a (real) degree 2 map $Y
\rightarrow \mathbb{P}^2$ whose branch locus is a nonsingular
quartic $C$.

\textbf{Case 2}: Here $Y_{\mathbb{R}}=2V_1$. According to
Theorem~\ref{p2}, $C_{\mathbb{R}}=\varnothing$ and
$Y_{\mathbb{R}}$ is the trivial double covering of
$\mathbb{P}_{\mathbb{R}}^2\cong V_1$. The deck translation
involution~$\Gamma$ of the covering permutes the two projective
planes. According to Lemma~\ref{lift2}, $\Gamma$ defines an
automorphism of an appropriate real Enriques surface and the
resulting automorphism realizes the permutation of the two copies
of $V_1$ in $E_{\mathbb{R}}^{(2)}$.

\textbf{Case 3}: In this case, $Y_{\mathbb{R}}=iS$, $i=2,3,4$.
According to Theorem~\ref{p2}, $C_{\mathbb{R}}$ consists of $i$
ovals for each $i=2,3,4$, and $Y_{\mathbb{R}}$ is the double
covering of the orientable part of
$\mathbb{P}_{\mathbb{R}}^2\setminus C_{\mathbb{R}}$ branched over
$C_{\mathbb{R}}$. According to Lemma~\ref{quartic}, any
permutation of the ovals of~$C$ can be realized by a rigid
isotopy. The latter defines a deformation of $Y$ and, via inverse
Donaldson's trick, a deformation of an appropriate real Enriques
surface that realizes the corresponding permutation of the spheres
of ${E_{\mathbb{R}}}^{(2)}$.

\textbf{Case 4}: The corresponding surface $Y$ obtained by
Donaldson's trick is a real unnodal Del Pezzo surface of degree 4
with  $Y_{\mathbb{R}}=2S$. According to the models of Del Pezzo
surfaces, the anti-canonical system $|-K|$ embeds $Y$ into
$\mathbb{P}^4$ as a surface of degree~$4$. Unnodal Del Pezzo
surfaces of degree~$4$ are regular intersections of two quadrics
in $\mathbb{P}^4$. From Theorem~\ref{p4}, all such surfaces in
$\mathbb{P}^4$ with real part $2S$ are rigidly isotopic. Thus, we
can take for~$Y$ the intersection of the $3$-sphere
$\{x^2+y^2+z^2+t^2=2u^2\}$ and the cylinder $\{x^2+y^2+z^2=u^2\}$,
where $[x:y:z:t:u]$ are homogeneous coordinates in $\mathbb{P}^4$.
Rotating the cylinder about the $tu$-plane through an angle $\pi$
realizes the desired permutation of the spheres.
\end{proof}

\begin{theorem} \label{almostDP}
With one exception, any permutation of homeomorphic components of
both halves of a real Enriques surface with a disconnected half
$E_{\mathbb{R}}^{(1)}= V_d \sqcup ...$, $d\geq 4$, can be realized
by deformations and automorphisms. In the exceptional case
$E_{\mathbb{R}}= \{V_4 \sqcup S \}\sqcup \{4S \}$, the realized
group is $D_{8}\subset S_4$.
\end{theorem}

\begin{remark} \label{no.others2}
In the exceptional case, the Pontrjagin-Viro form is well defined.
The quarter decomposition of~$E_\mathbb{R}^{(2)}$ is $(2S)\sqcup
(2S)$. A permutation of the spherical components is not realizable
if it doesn't preserve the quarter decomposition.
Theorem~\ref{almostDP} states that a permutation of the components
can be realized if and only if it preserves the quarter
decomposition.
\end{remark}

\begin{proof}
The problem reduces to a question about appropriate
$(g,r)$-surfaces, $g\geq 3$ and $r\geq 1$ (see~\cite {D2} for the
models of $(g,r)$-surfaces). We construct a particular surface
(within each deformation class) that has a desired automorphism or
`auto-deformation'. Among the extremal types listed in~\cite{D2},
we need to consider only the following types (as in the other
cases there are no homeomorphic components) and all their
derivatives $(E_{\mathbb{R}}^{(1)},\cdot )$ obtained from the
extremal ones by sequences of topological Morse simplifications of
$E_\mathbb{R}^{(2)}$:

\begin{enumerate}
\item $E_{\mathbb{R}}^{(1)}=V_4 \sqcup 2V_1$;\quad
$E_{\mathbb{R}}^{(2)}=\varnothing$; \item
$E_{\mathbb{R}}^{(1)}=V_4 \sqcup S$;\quad \quad
$E_{\mathbb{R}}^{(2)}=4S$.
\end{enumerate}

\textbf{Case 1}: In this case, homeomorphic components are in
$E_{\mathbb{R}}^{(1)}= B/t^{(2)}$ so we need to deal with $B$. By
Donaldson's trick, we obtain a \textit{DPN}-pair $(Y,B)$, where
$Y$ is a real $(3,2)$-surface with empty real part and $B\cong S_3
\sqcup 2S$ is an admissible branch curve on $Y$ such that the
rational components of $B$ are real. According to the models of
$(3,2)$-surfaces, $Y$ blows down to
$\Sigma_0={\mathbb{P}}^{1}\times {\mathbb{P}}^{1}$ with the real
structure $c_0\times c_1$, where $c_0$ is the usual complex
conjugation and $c_1$ is the quaternionic real structure ($\Fix
c_1=\varnothing$) on ${\mathbb{P}}^{1}$. The image $Q$ of $B$ is
the transversal union of smooth components $C'$, $C''$ and $C'''$,
where $C',C''\in |l_2|$ are two distinct real generatrices and
$C'''\in |4l_1+2l_2|$. By Theorem 18.4.1 in~\cite{D2}, there is
only one rigid homotopy class of such curves $Q\subset \Sigma_0$
and if $Q'$ is rigidly homotopic to $Q$ then the
\textit{DPN}-resolutions of the pairs $(\Sigma_0,Q')$ and
$(\Sigma_0,Q)$ are deformation equivalent in the class of
admissible \textit{DPN}-pairs. By a rigid homotopy of real
algebraic curves on $\Sigma_0$, we mean a path $Q_s$ of real
curves on $\Sigma_0$ such that the members of the path consist of
a fixed number of smooth components and have at most type $A$
singular points. Thus, a generic rigid homotopy of $Q_s$ defines a
deformation of the \textit{DPN}-resolutions $(Y_s,B_s)$ of the
pairs $(\Sigma_0,Q_s)$, thus, a deformation of the covering
$K3$-surfaces. Choosing a continuous family of invariant K\"ahler
metrics leads to a deformation of the corresponding real Enriques
surfaces obtained from inverse Donaldson's trick.

It suffices to connect $Q$ with itself by a path that realizes the
permutation of $C'$ and $C''$, and such that the members $Q_s$ of
the path split into sums $C_s'+C_s''+C_s'''$ of distinct real
smooth irreducible curves such that $C_s',C_s''\in |l_2|$ and
$C_s'''\in |4l_1+2l_2|$. Identify the real part of the base
$\mathbb{P}_{\mathbb{R}}^1\cong S^1=\mathbb{R}/2\pi$. Let
$Q=A_0+A_{\pi}+A$, where $A_\alpha$ is the generatrix of
$\Sigma_0$ over $\alpha$. Then, the family
$\{Q_t\}=\{A_t+A_{\pi+t}+A; t\in [0,\pi]\}$, defines a path that
realizes the permutation of the generatrices $A_0$ and $A_{\pi}$.\\

\textbf{Case 2}: Here we consider the 3 subcases:
$$E_{\mathbb{R}}=\{V_4 \sqcup S\} \sqcup \{iS\},~~i=2,3,4.$$
The corresponding \textit{DPN}-pair resulting from Donaldson's
trick is $(Y,B)$, where $Y$ is a real $(3,2)$-surface with
$Y_\mathbb{R}=iS$, and $B\cong S_3\sqcup 2S$ is an admissible
branch curve such that rational components of $B$ are conjugate
and $[B]=0$ in $H_2(X)$, where $X$ is the covering
\textit{K3}-surface. According to the models of $(3,2)$-surfaces,
there is a real regular degree $2$ map $\phi:Y\rightarrow
Z=\Sigma_4$ (with standard real structure, {\em i.e.},
$Z_{\mathbb{R}}=S_1$) branched over a nonsingular real curve $U\in
|2e_\infty|$. The irrational component of $B$ is mapped to a real
curve $F\in|e_\infty|$ and each rational component is mapped
isomorphically to the exceptional section $E_0$ of $Z$.
$U_{\mathbb{R}}$ is contained in a connected component of
$Z_{\mathbb{R}}\backslash((E_0)_\mathbb{R}\sqcup F_{\mathbb{R}})$.
By Theorem~\ref{Sig4}, up to rigid isotopy and automorphism of
$Z$, the pair $(U,F)$ is determined by its real root scheme. By
Theorem 18.4.2 in~\cite{D2}, the real \textit{DPN}-double of
$(Z;U,E_0\cup F)$ is determined up to deformation in the class of
admissible \textit{DPN}-pairs by the real root scheme of the pair
$(U,F)$. Thus, a generic rigid isotopy of the pairs $(U_s,F_s)$
defines a deformation of the \textit{DPN}-doubles $(Y_s,B_s)$ of
$(Z;U_s,E_0\cup F_s)$, so a deformation of the covering
$K3$-surfaces. Choosing a continuous family of invariant K\"ahler
metrics gives a deformation of the corresponding real Enriques
surfaces obtained from inverse Donaldson's trick. Thus, in view of
the above observation and Lemma~\ref{lift3}, it suffices to
realize permutations of certain ovals of $U$ by rigid isotopies of
the pair $(U,F)$ and/or involutive automorphisms of $Z$ preserving
$U$. In the latter case the set of fixed points should have
nonempty intersection with $U$.

The real root scheme of $(U,F)$ is a disjoint union of $i$
segments ({\em cf.} the first row of Table~\ref{root} for $i=4$),
and it has a representative (real root marking) with the desired
symmetry group ({\em i.e.}, $S_2$, $S_3$ and $D_8$ for $i=2,3$ and
4 respectively), generated by rotations and reflections of
$B_{\mathbb{R}}\cong S^1$. By Remark~\ref{common}, these
symmetries realize permutation of the corresponding ovals.
Furthermore, the fixed point set of a reflection symmetry consists
of two distinct points on $S^1$, which correspond to two distinct
real generatrices of $Z$. Since $U\in |2e_{\infty}|$, $U$
intersects the set of the fixed points of the induced involution
and one can apply Lemma~\ref{lift3}. For $i=4$, the reason, why
other permutations are not allowed is explained in
Remark~\ref{no.others2}.
\end{proof}

\begin{theorem} \label{(2,r)}
For the real Enriques surfaces with disconnected
$E_{\mathbb{R}}^{(1)}=V_3\sqcup ...$, none of the permutations of
the components of the half $E_{\mathbb{R}}^{(1)}$ is realizable by
deformations or automorphisms. With the exceptions listed below,
any permutation of homeomorphic components of the half
$E_{\mathbb{R}}^{(2)}$ can be realized by deformations and
automorphisms. The exceptional cases are:
\begin{enumerate}
\item surfaces with $E_\mathbb{R}=\{V_3\sqcup
V_{1}\}\sqcup\{4S\}:$ the realized group is $D_8;$ \item surfaces
with $E_\mathbb{R}=\{V_3\sqcup S\}\sqcup\{V_1\sqcup 3S\}:$ the
realized group is $S_2;$ \item surfaces with
$E_\mathbb{R}=\{V_3\sqcup V_1\sqcup S\}\sqcup\{3S\}:$ the realized
group is $S_2;$ \item surfaces with $E_\mathbb{R}=\{V_3\sqcup
2S\}\sqcup\{V_1\sqcup 2S\}:$ the realized group is trivial.
\end{enumerate}
\end{theorem}
\begin{remark}\label{no.others3}
In the exceptional cases, the Pontrjagin-Viro form is well
defined. It defines the quarter decompositions as follows:
\begin{enumerate}
\item $E_\mathbb{R}=\{(V_3\sqcup
V_{1})\sqcup(\varnothing)\}\sqcup\{(2S)\sqcup(2S)\}$; \item
$E_\mathbb{R}=\{(V_3\sqcup
S)\sqcup(\varnothing)\}\sqcup\{(V_1\sqcup S)\sqcup(2S)\}$; \item
$E_\mathbb{R}=\{(V_3\sqcup
S)\sqcup(V_1)\}\sqcup\{(2S)\sqcup(S)\}$; \item
$E_\mathbb{R}=\{(V_3\sqcup S)\sqcup(S)\} \sqcup\{(V_1 \sqcup
S)\sqcup(S)\}$.
\end{enumerate}
\noindent One cannot permute homeomorphic components without
preserving the quarter decomposition. The above theorem states
that a permutation of homeomorphic components of
$E_\mathbb{R}^{(2)}$ can be realized if and only if it preserves
the quarter decomposition.
\end{remark}

\begin{proof}
For these surfaces, the \textit{DPN}-pair resulting from
Donaldson's trick is $(Y,B)$, where $Y$ is a real unnodal
$(2,r)$-surface and $B\cong S_2\sqcup rS$ is an admissible branch
curve on $Y$, $r\geq 1$. According to the models of
$(2,r)$-surfaces (\cite{D2}), the anti-bicanonical system of $Y$
defines a surjective degree $2$ map $\varphi :Y \rightarrow
Z'\subset \mathbb{P}^3$. The branch locus of $\varphi$ is a cubic
section through the vertex. The map $\varphi$ lifts to a map
$\tilde{\varphi} :Y \rightarrow Z=\Sigma_{2}\subset
\mathbb{P}^3\times \mathbb{P}^1$ (with standard real structure)
and the branch locus of $\tilde{\varphi}$ is a curve $U \in
|2e_{\infty}+2l|$. By $\tilde{\varphi}$, the genus $2$ component
of $B$ is mapped to a nonsingular real generic section $F\in
|e_{\infty}|$ and the rational components of $B$ are mapped to the
exceptional section $E_0$ in $Z$. The pair $(U,F)$ is a suitable
pair (see Section~\ref{suip}). The pullback ${\tilde{\varphi}}^{-1}(E_0)$
consists of the fixed components of $|{-2K}|$ ($(-4)$-curves on
$Y$, {\em i.e.}, the rational components of $B$) and, possibly,
several $(-1)$-curves. Figure~\ref{dynkin} shows the Dynkin graph
of this configuration of curves. It is a linear tree with $2r-1$
vertices, where the two outermost ones, marked with $\tilde{E}_0$,
represent the components of the proper transform of $E_0$.

\begin{figure}[h]
\centering\epsfig{file=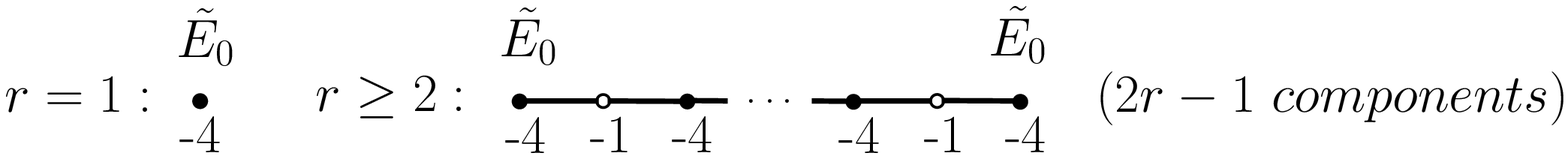,width=4.5in}\caption{}\label{dynkin}
\end{figure}

We start by proving the first part of the theorem. In view of
Theorem 19.1 in~\cite{D2}, we only need to consider the real
Enriques surfaces with $E_\mathbb{R}^{(1)}=V_3\sqcup mV_1\sqcup
nS$, $m=0$ or $1$ and $n=2, 3$ or $4$. By Donaldson's trick we
obtain real $(2,r)$-surfaces with admissible branch curves $B\cong
S_2\sqcup rS$, where $r=m+2n$, (as
$E_\mathbb{R}^{(1)}=B/t^{(2)}$). The Dynkin graph of the pullback
of the exceptional section $E_0$ contains $m+2n$ copies of
$(-4)$-curves that correspond to the spherical components of $B$.
Since the map $\varphi$ is anti-bicanonical, our model is
canonical and both the Dynkin graph and the corresponding Coxeter
diagram on the covering $K3$-surface are rigid. The only map that
can realize a permutation of the spherical components of $B$ is
the deck translation of the covering $\varphi$ which changes the
order of the curves in the Dynkin graph. But since the spherical
components permuted by the deck translation are identified by the
map $t^{(2)}$ on the covering $K3$-surface and
$E_\mathbb{R}^{(1)}=B/t^{(2)}$, the result follows.

Proof of the second part is based on suitable pairs.
Theorem~\ref{SP} states that, up to rigid isotopy and automorphism
of $Z$, a suitable pair $(U,F)$ is determined by its real root
scheme. In view of Theorem~\ref{gensp} and Lemma~\ref{lift4}, it
is enough to realize the permutations of certain ovals by rigid
isotopies of the pair $(U,F)$ and/or involutive automorphisms of
$Z$ preserving $U$, where in the latter case the set of fixed
points should intersect $U$. Proof is very similar to that of
Theorem~\ref{almostDP}, case 2. Among the extremal types listed
in~\cite{D2}, we need to consider only the following types and all
their derivatives $(E_{\mathbb{R}}^{(1)},\cdot )$ obtained from
the extremal ones by sequences of topological Morse
simplifications of $E_\mathbb{R}^{(2)}$:

\begin{enumerate}
\item $E_{\mathbb{R}}^{(1)}=V_3\sqcup V_1$;\quad \quad\quad
$E_{\mathbb{R}}^{(2)}=4S$; \item $E_{\mathbb{R}}^{(1)}=V_3\sqcup
S$;\quad \quad \quad \ $E_{\mathbb{R}}^{(2)}=V_1\sqcup 3S$; \item
$E_{\mathbb{R}}^{(1)}=V_3\sqcup V_1\sqcup S$;\quad \
$E_{\mathbb{R}}^{(2)}=3S$; \item $E_{\mathbb{R}}^{(1)}=V_3\sqcup
2S$;\quad \quad \quad $E_{\mathbb{R}}^{(2)}=V_1\sqcup 2S$; \item
$E_{\mathbb{R}}^{(1)}=V_3\sqcup V_1\sqcup 2S$;\quad
$E_{\mathbb{R}}^{(2)}=2S$.
\end{enumerate}
\vspace{0.3cm}

The extremal real root schemes of the pairs $(U,F)$ for these
cases are listed in Table~\ref{root}; the others are obtained by
removing several segments not containing a $\circ$~-vertex. As in
the proof of Theorem~\ref{almostDP}, one can realize each root
scheme by a sufficiently symmetric representative, constructing
the desired groups of permutations of the ovals of $U$,
see~Remark~\ref{common}. In the case of automorphisms, induced by
reflection symmetries, we observe that the fixed point set
consists of a pair of generatrices and intersects $U\in
|2e_{\infty}+2l|$; hence, Lemma~\ref{lift4} applies.
Remark~\ref{no.others3} explains why other permutations are not
realizable.

\end{proof}

\end{document}